
\documentstyle{amsppt}
\addto\tenpoint{\normalbaselineskip20pt\normalbaselines}
\hsize=6.0 true in
\vsize = 8.5 true in
\NoRunningHeads
\def\rn{Radon-Nikod\'ym}
\def\tuxy{T:U_X(\Sigma)\to Y}
\def\tuy{T:U(\Sigma)\to Y}
\def\mxy{m:\Sigma \to L(X,Y)}
\def\ifdm{\int f\,dm}
\def\ux{U_X(\Sigma)}
\def\mt{\widetilde m}
\def\tcxy{T:C(H,X)\to Y}
\def\tcy{T:C(H)\to Y}
\def\bh{\beta (H)}
\def\chx{C(H,X)}
\def\tp{T_\pi}
\def\tsc{T\sp \sharp:C(H)\to L(X,Y)}
\def\ran{\rangle}
\def\lan{\langle}
\def\cor{\leftrightarrow}
\def\sap{\sum_{A\in\pi}}
\def\mp{m_\pi}
\def\tshp{T\sp \sharp}
\def\mshp{m\sp \sharp}
\def\osm{(\Omega, \Sigma, \mu)}

\topmatter
 \title On Convergence of Conditional Expectation Operators \endtitle
 \author C. Bryan Dawson \endauthor
 \abstract{Given an operator $\tuxy$ or $\tcxy$, one may consider the
 net of conditional expectation operators $(\tp)$ directed by refinement
 of the partitions $\pi$.  It has been shown previously that $(\tp)$ does not
always converge to $T$. This paper gives several conditions under which this
convergence does occur, including complete characterizations when $X={\bold R}$
or when $X\sp *$ has the \rn\ property. }\endabstract
\endtopmatter
\document

 \specialhead 1\. Introduction \endspecialhead

 It is well known that if $\tuxy$ is a bounded linear operator, where
 $\ux$ is the uniform closure of the $X$-valued $\Sigma$-simple
 functions, then there is a unique finitely additive set function $\mxy$
 with finite semivariation such that $T(f)=\ifdm$ for all $f\in\ux$.
 Also, if $\tcxy$, there is a unique weakly regular
 $m:\bh\to L(X,Y\sp {**})$ such that $T(f)=\ifdm$ for $f\in\chx$.
 In each case, $\mt(H)=\|T\|$. Given such an operator $T$, a
 finite partition $\pi$ of $H$, and a measure $\mu$ on $H$, a conditional
 expectation operator $\tp$ can be defined. It was shown in \lbrack 1\rbrack
that
 the net $(\tp)$ directed by refinement does not always converge to $T$
 in the operator norm. Conditions under which this convergence does occur
are discussed herein.

 Throughout, $X$ and $Y$ are Banach spaces. The closed unit ball of $X$
 is denoted by $B_X$.  We will use $H$ for a compact Hausdorff space, and
$C(H,X)$ for the space
 of continuous functions from $H$ to $X$. An arbitrary $\sigma$-algebra
 of subsets of some universal space $\Omega$ will be represented by
 $\Sigma$, and when $\Omega=H$ we will use $\Sigma=\bh$, the Borel sets
 of $H$, without further mention. An additive set function $m:\Sigma\to
 X$ will be called a vector measure, while by a measure we mean a
countably additive set function
 $\mu:\Sigma\to\lbrack 0,\infty)$.

 For a vector measure $m:\Sigma \to X$, we define the variation of $m$ as
 usual and the scalar semivariation $\|m\|$ of $m$ as in \lbrack 8\rbrack .
 If $\mxy$, the semivariation $\mt$ of $m$ is as in \lbrack 10\rbrack  (recall
  that if $X={\bold R}$, then $\mt(A)=\|m\|(A)$,
 and if $Y={\bold R}$, then $\mt(A)=|m|(A)$).
We write $m\cor T$ to mean $m$ corresponds to $T$ as in the two
 theorems stated at the beginning of this section.

 If $m:\Sigma \to X$ is a
 vector measure, $\mu$ a positive, bounded, finitely additive measure on
 $\Sigma$, and $\pi$ a partition of $\Omega$, then the conditional
 expectation $\mp$ of $m$ by $\pi$ and $\mu$
is given by
   $$\mp(B)=\sap\frac{\mu(A\cap B)}{\mu(A)}m(A)\text{, observing }{\tsize\frac
00}=0.$$ If $m\cor\tuxy$ or
 $m\cor \tcxy$ is a bounded linear operator, $\pi$ is a partition, and
 $\mu$ is a measure on $\Sigma$ ($\beta(H)$), then the conditional expectation
operator $\tp$ is given by
  $$\tp(f)=\sap m(A)\left(\frac{\int_A
f\,d\mu}{\mu(A)}\right)\text{, again observing }{\tsize\frac 00}=0.$$
If $m\cor\tuxy$, then
 $\mp\cor\tp$ for each partition $\pi$. Also, if $m\cor\tcxy$ and
 $\mu$ is a regular measure on $\beta(H)$, then $\mp\cor\tp$; see \lbrack 5,
 Lemma 2.2\rbrack . We reserve $\tp\to T$ to mean convergence in operator
 norm. For other terminology not defined here, see \lbrack 8\rbrack  or
\lbrack 4\rbrack .

\specialhead 2\. The Results \endspecialhead

We now turn our attention to convergence of $\tp$ to $T$ in operator norm.
Lemma  1 gives a necessary condition for that convergence.

\proclaim{Lemma 1} Let $m\cor\tcxy$ (respectively $m\cor\tuxy$) and
let $\mu$ be a regular measure (a measure) on $\bh$ (\,$\Sigma$). If
$\tp\to T$, then $m$ has relatively compact range.\endproclaim

\demo{Proof} We have $\mp\to m$ in semivariation norm since $\tp\cor\mp$,
using the regularity of
$\mu$ in the $\chx$ case. Hence, $\mp(A)\to m(A)$
uniformly in $A$.  As the range of each $\mp$
is bounded and finite dimensional, we must have that the range of $m$ is
totally bounded.\enddemo

The converse is, in general, not true. A proof appears later. In the
examples of \lbrack 1\rbrack , the measures (with possibly one
exception) do not have relatively norm compact ranges. This lemma leads to the
following when $X={\bold R}$:

\proclaim{Theorem 2} Let $m\cor\tcy$ and suppose $m\ll\mu$. The
following are equivalent:
\roster\item $\mp\to m$ in semivariation norm.
\item $m$ has relatively compact range.
\item $T$ is compact.\endroster
If, in addition, $\mu$ is regular, then (1) -- (3) are
equivalent to
\roster\item\lbrack 4\rbrack  $\tp\to T$.\endroster
\endproclaim

\demo{Proof} The fact that (2) implies (1) in this setting follows from
\lbrack 3,
remark 5.2\rbrack .  By Lemma 1, (1) and (2) are equivalent. The equivalence
of (2)
and (3) is in \lbrack 11, p\. 496\rbrack . The proposition follows. \enddemo

The conditions (1)--(4) are also equivalent in the setting $m\cor\tuy$.
This theorem is applied frequently in the remainder of the paper.

Conditions under which $\tp\to T$ when $T$ takes its values in ${\bold R}$
are given after the following lemma. A martingale convergence theorem similar
to
the lemma can be found in \lbrack 7\rbrack .

\proclaim{Lemma 3} The following are equivalent:
\roster\item X has the \rn\ property (RNP).
\item If $\osm$ is a finite measure space,
$m:\Sigma\to X$ is of bounded variation and $m\ll\mu$, then
$\mp\to m$ in variation norm.
\endroster\endproclaim

\demo{Proof} That (1) implies (2) is seen from \lbrack 2,
Thm\. 1\rbrack  (and its
proof) applied to the singleton $\{m\}$. Now, suppose (2) holds and fix the
finite measure space $(\Omega,
\Sigma, \mu)$, and suppose $m:\Sigma\to X$ is of bounded variation and
$m\ll\mu$. For each partition $\pi$, define $f_\pi:\Omega\to X$ by
  $$f_\pi = \sap\frac{m(A)}{\mu(A)}\chi_A.$$
Then,
  $\mp(B)= \int_B
  f_\pi\,d\mu$ for all $B$, and
$|\mp|(\Omega)=\|f_\pi\|_1$ for each $\pi$. However, $(\mp)$
converges in variation norm by hypothesis; consequently $(f_\pi)$ is
Cauchy  and converges to some $f\in
L\sp 1(\mu,X)$. Then for each $B\in\Sigma$, $\int_B f_\pi\,d\mu\to
m(B)$, but also $\int_B f_\pi\,d\mu\to\int_B f\,d\mu$. Hence $m(B)=\int_B
f\,d\mu$ for all $B\in\Sigma$.   \enddemo

Note that Lemma 3 provides a converse of \lbrack 2, Thm\. 1\rbrack .
Consider applying the above theorem when $m\cor\tcxy$ or $m\cor\tuxy$;
conditions under which $\mp\to m$ in variation norm are easily
extracted. Also, if $Y=\text{{\bf R}}$ we have $\|T\|=\mt(\Omega)=|m|(\Omega)$,
which  proves the following theorem.

\proclaim{Theorem 4} The following are equivalent:
\roster\item $X\sp *$ has RNP.
\item $\tp\to T$ whenever $m\cor T:\chx\to\text{{\bf R}}$ and $m\ll\mu$.
\item $\tp\to T$ whenever $m\cor T:\ux\to\text{{\bf R}}$ and $m\ll\mu$.
\endroster\endproclaim

We shall now explore the general question for spaces whose duals possess RNP.

\proclaim{Definition} Let $m\cor\tuxy$ or
$m\cor\tcxy$. For $y\sp *\in Y\sp *$, define $m_{y\sp *}:\Sigma\to X\sp *$ by
$m_{y\sp *}(A)(x) = \lan y\sp *,m(A)(x)\ran$.\endproclaim

Note that $|m_{y\sp *}|(\Omega)=\widetilde{m_{y\sp *}}(\Omega)\le
\|y\sp *\|\mt(\Omega)$.

\proclaim{Theorem 5} The following are equivalent:
\roster\item $X\sp *$ has RNP.
\item $\tp\to T$ whenever $m\cor\tuxy$, where $T$ is compact and
$\{|m_{y\sp *}|:y\sp *\in B_{Y\sp *}\}$ is uniformly absolutely continuous with
respect to a positive, bounded, finitely additive $\mu$.
\item $\tp\to T$ whenever $m\cor\tcxy$ where $T$ is compact, $\mu$ is as in
(2), and $\mu$ is regular.\endroster\endproclaim

\demo{Proof} Suppose (1) holds, let $X, m$, and $T$ be as stated, and consider
$T\sp *:Y\sp *\to U_X(\Sigma)\sp *$. Now $U_X(\Sigma)\sp * =
 \{\mu\ |\ \mu:\Sigma\to
X\sp * \text{ has finite variation}\}$.
Also, $T\sp *(y\sp *)\cor m_{y\sp *}$. Since $T\sp *$ is compact,
 $\{m_{y\sp *}:y\sp *\in B_{Y\sp *}\}$ is relatively compact in variation
 norm. Let $\mu$
be as stated (existence is guaranteed by \lbrack 2\rbrack ). Applying \lbrack
2\rbrack  again,
 $(m_{y\sp *})_\pi\to
m_{y\sp *}$ uniformly in $y\sp *\in B_{Y\sp *}$, in variation norm. Hence
$T\sp *(y\sp *)_\pi\to T\sp *(y\sp *)$ uniformly in $y\sp *\in B_{Y\sp *}$.
 However, for
$f\in\ux$, we have
  $\tp\sp *(y\sp *)(f) =
  T\sp *(y\sp *)_\pi(f)$.
Thus, $\tp\sp *(y\sp *)\to T\sp *(y\sp *)$
uniformly in $y\sp *\in B_{Y\sp *}$, i.e\. $\tp\sp *\to T\sp *$. Therefore,
(2)
holds.

For (2) implies (3), consider $\widehat T:U_X(\bh)\to Y$. As $T$ is compact,
$m$ takes its values in $L(X,Y)$ and $\widehat T$ is also compact. By
hypothesis, $\widehat
{\tp}\to \widehat T$. As $\tp = \widehat {\tp}|_{\chx}$ and $T=\widehat
T|_{\chx}$, we have $\tp\to T$.

It remains to show that (3) implies (1). Let $m\cor T : C(\lbrack 0,1\rbrack
,X)\to
{\bold R}$ such that
$m\ll\lambda = $ Lebesgue measure. Let $y\in Y$ with $\|y\|=1$ and
define $T':C(\lbrack 0,1\rbrack ,X)\to Y$ by $T'(f)=T(f)y$. Define $m':\Sigma
\to
L(X,Y)$ by $m'(A)(x) = m(A)(x)y$, for all $A\in\Sigma$, $x\in X$. Then,
 $m'\cor T'$  and $|m'|\ll
\lambda$.  If $\|y\sp *\|\le 1$, then
$|m_{y\sp *}'|(A)\le|m'|(A)$, and consequently $\{|m_{y\sp *}| : y\sp *\in
B_{Y\sp *}\}$ is uniformly absolutely continuous with respect to $\lambda$.
Hence, by hypothesis, $\tp'\to T'$. But,
it can readily be seen that $\tp'(f)=\tp(f)y$; thus $\tp\to T$. By
Theorem 4, $X\sp *$ has RNP.\enddemo

 For a related theorem, see \lbrack 6, Theorem 2.3\rbrack .
Recall now the following
definition:

\proclaim{Definition} Suppose $m\cor\tcxy$. Define
$\tsc$ by $\tshp(f)(x) = T(xf)$ for $f\in
C(H)$. Define $\mshp:\Sigma\to L(\text{{\bf R}},L(X,Y))$ by $\mshp(A)(r) =
rm(A)$ for $r\in\text{{\bf R}}$, $A\in\Sigma$.\endproclaim

Note that if $m\cor\tcxy$ and $m$ is strongly bounded,
then $\mshp\cor\tshp$.

\proclaim{Corollary 6} Suppose $m\cor\tcxy$, where $X\sp *$ has RNP@. If
$T$ is compact, then so is $\tshp$.\endproclaim

\demo{Proof} Since $T$ is compact, it is strongly bounded (see \lbrack 4,
Thm\. 4.2\rbrack ). Then there is a regular measure $\mu$ on $\Sigma$ such that
$\mt\ll\mu$. By Theorem 5, $\tp\to T$.
Consequently $m$ (and hence $\mshp$) has relatively compact range by
Lemma 1. Finally, an application of Theorem 2 yields the
compactness of $\tshp$.\enddemo

In \lbrack 13\rbrack , Saab and Smith showed that $T$ nuclear implies $\tshp$
nuclear if and
only if $X\sp *$ has RNP@. It is not known if the converse of Corollary 6 is
true. However, we do have the following example:

\proclaim{Example 7} $T$ is compact need not imply $\tshp$ is compact.
\endproclaim

\demo{Demonstration} Let $m:\beta(\lbrack 0,1\rbrack )\to L(C(\lbrack 0,1
\rbrack ),\text{{\bf R}}) = C(\lbrack 0,1\rbrack )\sp *$ be
given by $m(A)(f) = \int_A f\,d\lambda$, where $\lambda$ is Lebesgue
measure. Then  $\|m(A)\|\le
\lambda(A)$, and therefore $m$ is a dominated representing measure (see
\lbrack 4, Thm\. 2.8\rbrack ). Let $T:C(\lbrack 0,1\rbrack ,C(\lbrack 0,1
\rbrack ))\to\text{{\bf R}}$ be given by
$T(f)=\ifdm$. Then $T$ is certainly compact; however, the range of $m$
is not relatively compact. To see this, consider the sets
$A_n=\cup_{i=1}\sp {2\sp {n-1}}C_{n,2i-1}$ where the $C_{n,i}$ are the dyadic
intervals. Let $r_n$ be the n$\sp {\text{th}}$ Rademacher function. Then
for $k<n$, we have
  $  \|m(A_n)-m(A_k)\| \ge \frac 12$.
By Theorem 2, $\tshp$ is not compact.\enddemo

The next proposition illustrates another use of Theorem 2. It is a special
case of \lbrack 12, Prop\. 3\rbrack .

\proclaim{Proposition 8} Suppose $L(X,Y)$ has the weak \rn\ property (wRNP) and
$m\cor\tcxy$
where $m$ is of bounded variation. Let $f$ be the Pettis-integrable \rn\
derivative of $m$, and let $f_\pi$ be as in Lemma 3, where
$\mu=|m|$. Then $f_\pi\to f$ in Pettis norm.\endproclaim

\demo{Proof} By the observation of Stegall on the range of an indefinite Pettis
 integral (see f.i.
\lbrack 9\rbrack ) and Theorem 2,  $\mp\sp \sharp\to\mshp$ in semivariation =
scalar
semivariation norm. Therefore, $\mp\to m$ in scalar semivariation norm. Thus,
$|x\sp *\mp - x\sp *m|(\Omega)\to 0$, uniformly in $x\sp *\in B_{L(X,Y)\sp *}$.
However, we have that $x\sp *f_\pi$ is the \rn\ derivative of $x\sp *\mp$ and
$x\sp *f$ is the \rn\ derivative of $x\sp *m$. Thus,
  $$\int_\Omega |x\sp *f_\pi - x\sp *f|\,dm = |x\sp *\mp - x\sp *m|(\Omega),$$
and $x\sp *f_\pi\to x\sp *f$ uniformly over $x\sp *\in B_{L(X,Y)\sp *}$.
Therefore,
$f_\pi\to f$ in Pettis norm.\enddemo

We now show that the converse of Lemma 1 is, in
general, false.

\proclaim{Proposition 9} If $m\cor T:C(H,X)\to\text{{\bf R}}$, then
$\tshp_\pi\to\tshp$ need not imply that $\tp\to T$.\endproclaim

\demo{Proof} Let $X$ be a Banach space such that $X\sp *$ has wRNP but not
RNP (for example, $X$ can be the James tree space;
see \lbrack 14, p\. 87\rbrack  and \lbrack 8, p\. 214\rbrack ). Let $H=\lbrack
 0,1\rbrack $. Again by Stegall's
observation,
$\tshp$ is compact, for all $T:C(H,X)\to\text{{\bf R}}$. Let $\mu$ be Lebesgue
measure and let $m\ll\mu$ where $m\cor T:C(H,X)\to\text{{\bf R}}$. We then have
$\tshp_\pi\to\tshp$. If this implied that $\tp\to T$, then by Theorem
4, we would have that $X\sp *$ has RNP, a contradiction.\enddemo

Now consider the first three statements of Theorem 2 in the setting
$m\cor\tcxy$ with $m\ll\mu$, $\mu$ regular. The only implication that
remains true in this setting is (1) implies (2), provided by Lemma 1.
The proof of Proposition 9 above shows that (2) and (3) together need
not imply (1), and the fact that (3) does not imply (2) was treated in
Example 7. To see that (1) and (2) does not imply (3), let $H$ be a
singleton and let $X$ and $Y$ be such that there is a noncompact
operator $f:X\to Y$. Let $m(H)=f$ and $m(\phi)=0$. Let $T\cor m$. Then
$C(H,X)\cong X$, and $\tcxy$ is actually $f:X\to Y$. Hence, $T$ is not
compact. However, $\mp=m$ for the unique partition $\pi$ of $H$ and $m$
certainly has relatively compact range.

\subhead Acknowledgement \endsubhead The author is grateful to his
advisor Professor Paul W. Lewis and to Professor Elizabeth Bator, both
of the University of North Texas, for many helpful discussions on this
material.

\addto\tenpoint{\normalbaselineskip10pt\normalbaselines}
\parindent=0pt
\parskip = 0 pt
Emporia State University

\vskip -10 pt Emporia, KS 66801

\enddocument